\documentclass[a4paper,11pt]{article}

\usepackage{amsmath, amscd, theorem}
\usepackage[numbers]{natbib}
\usepackage{verbatim}

\usepackage{amssymb}
\newcommand{\wtF}{\widetilde F}

\newcommand{\wtf}{\widetilde f}
\newcommand{\wtg}{\widetilde g}
\newcommand{\bra}{\langle }
\newcommand{\ket}{\rangle }

\newcommand{\cL}{\mathcal{L}}

\newcommand{\sfg}{\mathsf{g}}
\newcommand{\wsg}{\widetilde\sfg}

\theorembodyfont{\normalfont\slshape}

{\theorembodyfont{\rmfamily}
}
\theoremstyle{break}

\usepackage[numbers]{natbib}
\title{Densities, Laplace Transforms and Analytic Number Theory}
\author{ Sibusiso Sibisi \thanks{\texttt{ssibisi@csir.co.za}} \\
      CSIR, Pretoria, South Africa \\
      \texttt{www.csir.co.za} }
\date{October 2008}

\begin{document}
\maketitle
\thispagestyle{empty}

\begin{abstract}
\noindent Li showed that the Riemann Hypothesis is equivalent to the nonnegativity of a 
certain sequence of numbers. 
\citeauthor{Bombieri} gave an arithmetic formula for the number sequence based on the Guinand-Weil explicit formula and showed that Li's criterion is equivalent to Weil's criterion for the Riemann Hypothesis.
We provide a derivation of the explicit formula based on Laplace transforms and present an alternative expression for Li's criterion that invites a probabilistic interpretation.
\end{abstract}

\begin{flushleft}
\end{flushleft}

\section{Introduction}
The Laplace transform of a function $f(t)$ on $[0,\infty)$ is defined by
\begin{equation*}
 \widetilde f(s) \equiv \cL\{f\}(s)= \int_0^\infty e^{-st} f(t) dt 
\label{eq:lplc}
\end{equation*}
If $\wtf$ converges for $\Re s=s_0$ then it converges for all $s$ with $\Re s > s_0$.  

\begin{table}[tbh]
\renewcommand{\arraystretch}{2.0}

\[
\begin{array}{|l|l|l|l|l|l|l|}


\cline{2-3} \cline{6-7}
\multicolumn{1}{c}{} & 
\multicolumn{2}{|c|}{\hbox{General Properties}}  &
\multicolumn{1}{c}{\phantom{xxx}} & 
\multicolumn{1}{c}{} & 
\multicolumn{2}{|c|}{\hbox{Specific Cases}} \\ 
\cline{1-3} \cline{5-7}

{\rm G1} & e^{-at}f(t) & \wtf(s+a) & &
{\rm S1} & \dfrac{t^n}{n!} \quad (n\ge 0)  & \dfrac{1}{s^{n+1}} \\

{\rm G2} & f'(t) & s\wtf(s) - f(+0) & &
{\rm S2} & e^{at} & \dfrac{1}{s-a} \\

{\rm G3} &  t^n f(t) & (-)^n \wtf^{(n)}(s) & & 
{\rm S3} & \delta(t-a) \; (a\ge0) & e^{-as} \\

{\rm G4} & \displaystyle \int_0^t f(u) du & \dfrac{1}{s} \wtf(s)  & &
{\rm S4} & \displaystyle \sum_{n=0}^\infty \delta(t-na) & \dfrac{1}{1-e^{-as}} \\

\cline{5-7}

{\rm G5} &  \dfrac{1}{a} f\left(\dfrac{t}{a}\right) & \wtf(as) &  
\multicolumn{1}{c}{} &
\multicolumn{3}{c}{\delta(t) = {\hbox{Dirac delta function}}} \\

{\rm G6} &  (f*g)(t) & \wtf(s)\,\wtg(s) &  
\multicolumn{4}{c}{} \\

\cline{1-3} 

\end{array}
\]

\renewcommand{\arraystretch}{1.0}

\caption{Laplace Transform pairs:- $\wtf(s) = \int_0^\infty e^{-st} f(t) dt$}
\label{tbl:lplc}
\end{table}  Table~\ref{tbl:lplc}, sourced (directly or indirectly) from the more comprehensive table in \cite[p. 1020]{AbrSteg}, provides Laplace transform pairs that are relevant to this paper.
$f*g$ denotes Laplace convolution
\begin{equation*}
 (f*g)(t) = \int_0^t f(u) g(t-u)du
\end{equation*}
If $f(t)$ is a real, non-negative function then $F(t)=\int_0^t f(u) du$ is called a {\it distribution function} with {\it density} $f$.
$F$ is a probability distribution and $f$ the corresponding probability density if $F(\infty) = \wtf(0) = 1$. 

We take the view that the density $f$ is the fundamental object of study because it allows the construction of  a variety of (weighted) integrals involving $f$ over prescribed intervals of $[0,\infty)$, 
the distribution $F(t)$ being one such instance. 
This provides a simple, general and  coherent  framework for constructing explicit formulae of analytic number theory.

The paper is largely expository in nature, using Laplace transforms to reproduce known results in analytic number theory with what might arguably be regarded as natural ease.
The approach taken here is possibly best illustrated by example.
Accordingly, we shall first study the Chebyshev counting function before proceeding to the general case.

\section{The Chebyshev function: arithmetic form}
For $\Re s>1$, the Riemann zeta function is given by \cite[p. 6]{Edwards}
\begin{equation}
 \zeta(s) = \sum_{n=1}^\infty \frac{1}{n^s} =  \prod_q \frac{1}{1-q^{-s}} 
\label{eq:zeta}
\end{equation}
The rightmost form is the Euler product representation, where $q$ runs over all the primes $\{2,3,5,7,11,\ldots\}$.
Its logarithmic derivative gives
\begin{equation}
 -\frac{\zeta'(s)}{\zeta(s)} = \sum_q \log q \, \frac{q^{-s}}{1-q^{-s}} =
  \sum_q \log q \left(\frac{1}{1-q^{-s}}-1\right) 
\label{eq:zetalogderiv1}
\end{equation}
With the aid of S4 in Table~\ref{tbl:lplc}, we recognise this as the Laplace transform of 
\begin{equation}
 f(t) = \sum_q \log q \sum_{n=1}^\infty \delta(t-n\log q) \qquad t\ge0
\label{eq:density1}
\end{equation}
Hence $f$ is a discrete density with an atom of strength $\log q$ at every multiple of $\log q$. 
The associated distribution $F$ is given by
\begin{equation}
 F(t) =  \int_0^t f(u) du = \sum_q \log q \sum_{n=1}^\infty \int_0^t \delta(u-n\log q) du 
      = \sum_{q^n\le e^t} \log q
\label{eq:chebF1}
\end{equation}
The rightmost sum is a shorthand for a double sum over both $q$ and $n$ such that $q^n\le e^t$.
Setting $x=e^t$, the Chebyshev counting function $\psi(x)$ is defined by
\begin{equation}
 \psi(x) \equiv F(\log x) = \sum_{q^n\le x} \log q
\label{eq:psi1}
\end{equation}
It is often written in terms of the von Mangoldt function $\Lambda(n)$ as
\begin{equation}
 \psi(x) = \sum_{n\le x}\Lambda(n) \qquad
\Lambda(n) =
\begin{cases}
 \log q & n=q^k \\
0 & {\rm otherwise}
\end{cases}
\label{eq:psi1a}
\end{equation}
In words, $\psi(x)$ is a distribution function over the integers with jumps of size $\log q$ at every power $q^k\le x\,(k\ge 1)$.
The density (\ref{eq:density1}) can correspondingly be written 
as
\begin{equation}
 f(t) = \sum_{n=1}^\infty \Lambda(n) \delta(t-\log n) \qquad t\ge0
\label{eq:density1a}
\end{equation}
For completeness, we note that (\ref{eq:zetalogderiv1}) can also be written as 
\begin{equation}
 -\frac{\zeta'(s)}{\zeta(s)} = \sum_{n=1}^\infty \frac{\Lambda(n)}{n^s}
\label{eq:zetalogderiv1a}
\end{equation}

The results presented here are well-known. 
The point to be noted is that we have arrived at $\psi$ in a bottom-up fashion via the density $f$.
We repeat the exercise below for the analytic representation of the zeta function. 

\section{The Chebyshev function: analytic form}
The Riemann zeta and $\xi$ functions are related by \cite[p. 16]{Edwards}:
\begin{equation}
2\xi(s) = s(s-1) \pi^{-s/2}\Gamma(s/2) \zeta(s)
\label{eq:xizeta}
\end{equation}
Both $\xi$ and $\Gamma$ functions are expressible as infinite products
 \begin{equation*}
2\xi(s) = \prod_\rho \left(1-\frac{s}{\rho}\right) \qquad{\rm and}\qquad
s e^{\gamma s}\, \Gamma(s) =  \prod_{n=1}^\infty \left(1+\frac{s}{n}\right)^{-1} e^{s/n}
\end{equation*}
where $\gamma$ is the Euler-Mascheroni constant and the product formula for $\xi$ runs over all zeros $\{\rho\}$ of $\xi$ (complex zeros of $\zeta$)  with $\rho$ and $1-\rho$ paired together. 
This formulation of the zeta function is valid for all $s$ except for a simple pole at $s=1$.

Following the steps of the previous section, the logarithmic derivative of (\ref{eq:xizeta}) gives
\begin{equation}
 -\frac{\zeta'(s)}{\zeta(s)} =  \frac{1}{s-1}
 -\sum_{n=1}^\infty\left(\frac{1}{s+2n} - \frac{1}{2n}\right) 
 - \sum_\rho\frac{1}{s-\rho}
 -\frac{1}{2}(\gamma+\log\pi)
\label{eq:zetalogderiv2}
\end{equation}
Laplace inversion leads to an alternative form for the density
\begin{equation}
f(t) = e^t 
 - \sum_{n=1}^\infty\left(e^{-2nt} - \frac{1}{2n}\,\delta(t)\right)
 -\sum_\rho e^{\rho t}
 -\frac{1}{2}(\gamma+\log\pi)\delta(t) \qquad t\ge0
\label{eq:density2}
\end{equation}
Hence the corresponding distribution $F$ is
\begin{alignat}{1}
F(t) =  \int_0^t f(u) du 
    &= e^t +\sum_{n=1}^\infty \frac{e^{-2nt}}{2n}
          -\sum_\rho \frac{e^{\rho t}}{\rho}
        -1 + \sum_\rho \frac{1}{\rho} -\frac{1}{2}(\gamma+\log\pi)   
\nonumber \\
    &= e^t +\sum_{n=1}^\infty \frac{e^{-2nt}}{2n}
          -\sum_\rho \frac{e^{\rho t}}{\rho}
        -\log 2\pi \qquad t>0 
\label{eq:chebF2a}
\end{alignat}
The second form follows from setting $s=0$ in (\ref{eq:zetalogderiv2}), so that
\begin{equation}
-\frac{\zeta'(0)}{\zeta(0)} = - 1 + \sum_\rho\frac{1}{\rho} - \frac{1}{2}(\gamma+\log\pi)
 = -\log 2\pi
\label{eq:zetalogderiv0}
\end{equation}
For a proof of the rightmost equality, see \cite[p. 67]{Edwards}. 

Finally we may write $\psi(x)$ as
\begin{equation}
 \psi(x) \equiv F(\log x) =  x
 + \sum_{n=1}^\infty \frac{x^{-2n}}{2n}
 -\sum_\rho \frac{x^{\rho}}{\rho}
 - \log 2\pi \qquad x>1
\label{eq:psi2}
\end{equation}
This is known as von Mangoldt's explicit formula for the Chebyshev function $\psi(x)$.
The top-down approach to deducing the density from the distribution $\psi$ or $F$ is by differentiation
\begin{alignat*}{2}
 \psi'(x)  &=  1 - \sum_{n=1}^\infty x^{-2n-1} - \sum_\rho x^{\rho-1} \qquad & &x>1 
\\
{\rm or} \qquad F'(t) &= e^t - \sum_{n=1}^\infty e^{-2nt} -\sum_\rho e^{\rho t} \quad & &t>0 
\end{alignat*}
While this may be valid for $t>0$, the density remains unspecified at $t=0$ and therefore technically incomplete.
But, as will emerge below, the atom at zero plays a significant role. 

To put the foregoing discussion in different words, we first note that, by G3 in Table~\ref{tbl:lplc}
\begin{equation*}
\wtf(s) = -\frac{\zeta'(s)}{\zeta(s)} \quad \implies \quad
\wtF(s) = \frac{1}{s}\wtf(s) = -\frac{1}{s}\frac{\zeta'(s)}{\zeta(s)}
\end{equation*}
Starting from $\wtF(s)$, we can
\begin{enumerate}
 \item 
Invert $\wtF(s)$\footnote{This is the Laplace transform 
equivalent of Edwards' explicit inversion \cite[p. 50]{Edwards} to obtain $\psi(x)$ directly.} 
to obtain $F(t)$ and then differentiate to obtain the density $F'(t)$ valid for $t>0$. 
\item
Invert $s\wtF(s)$ to obtain the density directly.
Here, we use G2 in Table~\ref{tbl:lplc}:
\begin{equation*}
\cL^{-1}\{s\wtF(s) - F(+0)\} = F'(t) \quad\implies\quad
\cL^{-1}\{s\wtF(s)\} =   F'(t)+ F(+0)\delta(t)
\end{equation*}
Hence the density we seek is $F'(t)+F(+0)\delta(t)$, which is defined for all $t\ge0$. 
This is equivalent to the density (\ref{eq:density2}), with appropriate grouping of terms.
\end{enumerate}
Armed with a fully specified density, we now turn to the derivation of a more general explicit formula,
subsuming von Mangoldt's explicit formula (\ref{eq:psi2}) as a special case.

\section{The general explicit formula}
Given a function $w(t)$ on $[0,\infty)$, $\int_0^T w(t)f(t)dt$ can be written in two ways using the arithmetic and analytic forms of the density $f(t)$ given in (\ref{eq:density1}) and (\ref{eq:density2}) respectively. 
We have already dealt in detail with the case $w(t)=1$ above, where the integral yields the Chebyshev counting fuction in arithmetic and analytic form. 

In this section we shall primarily be interested in the limiting case  $T\rightarrow\infty$.
Subject to convergence, let $\bra w(t)\ket$
 denote the expectation of $w(t)$ with respect to $f(t)$
\begin{equation}
 \bra w(t)\ket = \lim_{T\rightarrow\infty}\int_0^T w(t) f(t) dt = \int_0^\infty w(t) f(t) dt 
\label{eq:xpctn}
\end{equation}
Consider $w(t) = e^{-st}g(t)$ where $g$ is a function on $[0,\infty)$ with Laplace transform $\wtg$.
The expectation with respect to the arithmetic density (\ref{eq:density1}) gives
\begin{alignat}{1}
 \langle e^{-st} g(t)\rangle  
   &= \lim_{x\rightarrow\infty}\int_0^{\log x} e^{-st} g(t) f(t) dt \qquad (T=\log x) \nonumber \\
   &= \lim_{x\rightarrow\infty}\sum_{n\le x} \Lambda(n) \int_0^\infty e^{-st} g(t) \delta(t-\log n)dt
\nonumber \\
   &= \lim_{x\rightarrow\infty}\sum_{n\le x} \frac{\Lambda(n)}{n^s} g(\log n)  
\label{eq:xplct1} 
\end{alignat}
The domain of validity in $s$ will depend on the choice of $g(t)$, 
{\it e.g.} $g(t)=e^{-\alpha t}$ requires $\Re s>1-\alpha$.
The expectation with respect to the analytic density (\ref{eq:density2}) gives
\begin{alignat}{1}
 &\langle e^{-st} g(t)\rangle  = \int_0^\infty e^{-st} g(t)f(t) dt \nonumber \\
 &= \int_0^\infty e^{-st} g(t) \left(e^t 
 - \sum_{n=1}^\infty\left(e^{-2nt} - \frac{1}{2n}\,\delta(t)\right)
 -\sum_\rho e^{\rho t}
 -\frac{1}{2}(\gamma+\log\pi)\delta(t)\right) dt  \nonumber \\
  &= \wtg(s-1) - \sum_{n=1}^\infty\left(\wtg(s+2n) - \frac{1}{2n}\,g(0)\right)
 -\sum_\rho \wtg(s-\rho)
 -\frac{1}{2}(\gamma+\log\pi) g(0)
\label{eq:xplct2}
\end{alignat}
This is the general explicit formula.
The domain of validity in $s$ will again depend on the choice of $g(t)$, although analytic continuation allows convergence in (\ref{eq:xplct2}) for regions of $s$ where (\ref{eq:xplct1}) does not converge.
We shall formally explore convergence for particular `test functions' in the sections below,
including the question of when we may legitimately equate (\ref{eq:xplct1}) and (\ref{eq:xplct2}). 
For ease of exposition,  we shall proceed with the general case for now 
without repeatedly calling to question its convergence properties.
We shall have a specific interest in the case $s=1$, so that (\ref{eq:xplct2}) becomes
\begin{equation}
  \sum_\rho \wtg(1-\rho) =
  \wtg(0) - \sum_{n=1}^\infty\left(\wtg(1+2n) - \frac{1}{2n}\,g(0)\right)
  -\frac{1}{2}(\gamma+\log\pi) g(0)
  - \langle e^{-t} g(t)\rangle 
\label{eq:xplct2s1}
\end{equation}

We proceed to derive an associated formula, using a related function $\sfg$ equivalent to what is referred to in \cite{Bombieri} as the involution of $g$.
We first note that G5 in Table~\ref{tbl:lplc} 
is normally defined for $a>0$.
It is also valid for $a<0$ provided that $f(t/a)$ and $\wtf(as)$ remain meaningfully defined.
For $a=-1$, define 
\begin{equation*}
 \sfg(t) = -e^tg(-t) \qquad\implies\qquad \wsg(s) = \wtg(1-s)
\end{equation*}
Hence
\begin{alignat}{1}
  &\langle e^{-st} \sfg(t)\rangle = -\int_0^\infty e^{-st} g(-t) f(t)dt \nonumber \\
  &= \wtg(2-s) - \sum_{n=1}^\infty\left(\wtg(1-s-2n) + \frac{1}{2n}\,g(0)\right)
 -\sum_\rho \wtg(1-s+\rho)
  + \frac{1}{2}(\gamma+\log\pi) g(0)
\label{eq:xplct2a}
\end{alignat}
Setting $s=1$ gives
\begin{equation}
  \sum_\rho \wtg(\rho) =
  \wtg(1) - \sum_{n=1}^\infty\left(\wtg(-2n) + \frac{1}{2n}\,g(0)\right)
  +\frac{1}{2}(\gamma+\log\pi) g(0)
  + \langle g(-t)\rangle
\label{eq:xplct2as1}
\end{equation}
Since $\sum_\rho \wtg(\rho) =  \sum_\rho \wtg(1-\rho)$ by virtue of the pairing of $\rho$ and $1-\rho$,
(\ref{eq:xplct2s1}) and (\ref{eq:xplct2as1}) are equivalent.

\section{Weil's criterion}
We now turn to the convolution $(\sfg * g)(t)$ whose Laplace transform is, by G6 in Table~\ref{tbl:lplc}, 
$\wsg(s)\wtg(s) = \wtg(1-s)\wtg(s)$.
Hence, noting that $(\sfg * g)(0)=0$,  the expectation of $\sfg * g$ is
\begin{alignat}{1}
 &\langle e^{-st} (\sfg *g)(t)\rangle  = \int_0^\infty e^{-st} (\sfg *g)(t)f(t) dt \nonumber \\
  &= \wtg(2-s)\wtg(s-1) - \sum_{n=1}^\infty \wtg(1-s-2n) \wtg(s+2n)
 -\sum_\rho \wtg(1-s+\rho) \wtg(s-\rho)
\label{eq:xplct2conv}
\end{alignat}
Setting $s=1$ as before gives
\begin{equation}
 \sum_\rho \wtg(\rho) \wtg(1-\rho)
  = \wtg(1)\wtg(0) - \sum_{n=1}^\infty \wtg(-2n) \wtg(1+2n) - \langle e^{-t} (\sfg *g)(t)\rangle
\label{eq:xplct2convs1}
\end{equation}
Weil's criterion states that a necessary and sufficient condition for the Riemann hypothesis to hold true is
\begin{equation}
 \sum_\rho \wtg(\rho) \wtg(1-\rho) > 0
\label{eq:weil}
\end{equation}
for all smooth functions $g(t)$ \cite{Bombieri}\footnote{The criterion generalises to complex smooth functions $g(t)$ but real $g(t)$ suffices for our purposes.}.

It is clearly not feasible to test the criterion for all conceivable smooth functions.
It is desirable therefore to find a manageable subset of test functions that can be shown to suffice.
Li's criterion,  discussed below, provides such a subset.
By way of motivation, consider functions $g$ that satisfy 
\begin{equation}
\begin{split}
(g*\sfg)(t) &= g(t)+\sfg(t)\; = g(t)-e^t g(-t) \\
\implies\quad \wtg(s)\wsg(s)  &= \wtg(s)+\wsg(s) = \wtg(s)+\wtg(1-s)
\end{split}
\label{eq:sum=prod}
\end{equation}
Then the left side of (\ref{eq:xplct2s1}) and (\ref{eq:xplct2as1}) can be written as 
\begin{equation}
 \sum_\rho \wtg(\rho) =  \sum_\rho \wtg(1-\rho) 
 = \frac{1}{2} \sum_\rho \wtg(\rho) + \wtg(1-\rho)
 = \frac{1}{2} \sum_\rho \wtg(\rho) \wtg(1-\rho)
\end{equation}
and (\ref{eq:xplct2convs1}) reduces to the sum of (\ref{eq:xplct2s1}) and (\ref{eq:xplct2as1}).
Hence, if a subclass of smooth functions satisfying (\ref{eq:sum=prod}) can be shown to suffice for  Weil's criterion, then the criterion amounts to the positivity of (\ref{eq:xplct2s1}) or, equivalently, (\ref{eq:xplct2as1}) for such a set of functions.

We now turn to some simple test functions to illustrate the ideas discussed thus far and to build toward the test functions needed for Li's criterion.

\section{Polynomial test functions}
The simplest case $g(t)=1$  reproduces $\langle e^{-st} g(t)\rangle = \wtf(s) = -\zeta'(s)/\zeta(s)$ which we already know can be written in both arithmetic (\ref{eq:zetalogderiv1a}) and analytic form 
(\ref{eq:zetalogderiv2}) for $\Re s > 1$
\begin{equation}
\sum_{n=1}^\infty \frac{\Lambda(n)}{n^s} - \frac{1}{s-1} 
 = - \sum_{n=1}^\infty\left(\frac{1}{s+2n} - \frac{1}{2n}\right)
 -\sum_\rho \frac{1}{s-\rho}
 -\frac{1}{2}(\gamma+\log\pi)
\label{eq:g=1}
\end{equation}
The choice $g(t)=t^{k},\; k>0$ gives $(-)^{k}\wtf^{(k)}(s)$
\begin{equation}
 \sum_{n=1}^\infty \frac{\Lambda(n)}{n^s} (\log n)^k - \frac{k!}{(s-1)^{k+1}} 
  = - k! \left(\sum_{n=1}^\infty\frac{1}{(s+2n)^{k+1}}
 +\sum_\rho \frac{1}{(s-\rho)^{k+1}}\right)
\label{eq:g=t^k}
\end{equation}
This shows that if $g(t)$ is a constant, a power of $t$ or a polynomial in $t$, then the explict formula amounts to working with derivatives of $\wtf(s)$ or a linear combination thereof for the polynomial case.
This may seem to suggest that we can dispense with the explicit formula.
However, the inadmissibility of $s=1$ in (\ref{eq:g=1}) and (\ref{eq:g=t^k}) --
where we have grouped together on the left side the terms that are not bounded at $s=1$ --
illustrates the benefit of retaining the view that the density $f$ is the fundamental object of study.

At $s=1$ we  retreat to the limit form (\ref{eq:xpctn}) for the left side so that (\ref{eq:g=t^k}), say, becomes
\begin{alignat}{1}
 &\lim_{x\rightarrow \infty} \left[\sum_{n\le x}\frac{\Lambda(n)}{n^s} (\log n)^k   
- \int_0^{\log x} t^k e^{-(s-1)t}dt \right]_{s=1} \nonumber \\
= &\lim_{x\rightarrow \infty} \left[\sum_{n\le x}\frac{\Lambda(n)}{n} (\log n)^k  
 - \int_0^{\log x} t^k dt \right] \nonumber \\
= &\lim_{x\rightarrow \infty} \left[\sum_{n\le x}\frac{\Lambda(n)}{n} (\log n)^k  
 - \frac{(\log x)^{k+1}}{k+1} \right] =  (-)^k k! \,\eta_k 
\label{eq:eta}
\end{alignat}
We have recognised the expression as (proportional to) the number $\eta_k$ that arises in the Laurent expansion\footnote
{We avoid power series in $s$, in keeping with the following sentiment that 
\citeauthor{Edwards} \cite[p. 9]{Edwards} attributes to Riemann:
``The view of analytic continuation in terms of chains of disks and power series convergent in each disk
descends from Weiersrass and is quite antithetical to Riemann's basic philosophy that analytic functions
should be dealt with {\em globally}, not locally in terms of power series.''}
about $s=1$ of $\zeta'(s)/\zeta(s)$ \cite{Bombieri}.
Hence (\ref{eq:g=t^k}) at $s=1$  is
\begin{equation}
(-)^k \,\eta_k  
= - \sum_{n=1}^\infty\frac{1}{(1+2n)^{k+1}} - \sum_\rho \frac{1}{(1-\rho)^{k+1}} \qquad k>0
\label{eq:g=t^ks1}
\end{equation}
It is straightforward to verify that
\begin{equation*}
\sum_{n=1}^\infty \frac{1}{(1+2n)^k} 
= \sum_{n=1}^\infty \left(\frac{1}{n^k} - \frac{1}{(2n)^k}\right) - 1
= \left(1-2^{-k}\right)\zeta(k) - 1 \qquad k>1
\end{equation*}
Hence (\ref{eq:g=t^ks1}) may be written as
\begin{alignat}{1}
\sum_\rho \left(-\frac{1}{\rho}\right)^{k+1} 
&= \eta_k + (-)^k (1-2^{-k-1})\zeta(k+1) - (-)^k
\label{eq:g=t^ks1a} \qquad k>0 \\
\intertext{The $k=0$ case is already known from (\ref{eq:zetalogderiv0}) to be}
\sum_\rho \frac{1}{\rho} 
&= 1 + \frac{1}{2}\left(\gamma -\log4\pi \right)
\label{eq:g=1s1} 
\end{alignat}
Finally, we consider $\sfg(t) = -e^tg(-t)$ so that, for $g(t) = t^k,\quad k\ge0$
\begin{equation*}
 \langle e^{-st} \sfg(t)\rangle = 
 - \wtf^{(k)}(s-1) = \frac{d^k}{ds^k}\frac{\zeta'}{\zeta}(s-1) 
\end{equation*}
Hence, for $k>0$, (\ref{eq:xplct2a}) becomes
\begin{equation}
\langle e^{-st} \sfg(t)\rangle 
= k!\left(\frac{1}{(2-s)^{k+1}} - \sum_{n=1}^\infty \frac{1}{(1-s-2n)^{k+1}}
  - \sum_\rho \frac{1}{(1-s+\rho)^{k+1}}\right) \qquad k>0
\label{eq:g=t^ka} 
\end{equation}
The $s=1$ case does not require a limiting approach this time since all terms are bounded.
Let 
\begin{equation}
 \mu_k = \frac{(-)^{k+1}}{k!}\langle e^{-t} \sfg(t)\rangle = \frac{(-)^k}{k!}\wtf^{(k)}(0) 
 = \frac{(-)^{k}}{k!}\frac{d^k}{ds^k}\left[-\frac{\zeta'(s)}{\zeta(s)}\right]_{s=0} \qquad k\ge0
\label{eq:mu}
\end{equation}
Then (\ref{eq:g=t^ka}) becomes
\begin{equation}
\sum_\rho \left(-\frac{1}{\rho}\right)^{k+1}
  = -\mu_k  - 2^{-(k+1)}\zeta(k+1) -  (-)^k \qquad k>0
\label{eq:g=t^kas1}
\end{equation}
The $k=0$ case simply reproduces (\ref{eq:g=1s1}).
For completeness, we note  that $\mu_0 = -\log\pi$.

\section{Li's Criterion}
Li \cite{Li} proved that a necessary and sufficient condition for the Riemann hypothesis to hold is
\begin{equation*}
 \lambda_n = \sum_\rho 1-\left(1-\frac{1}{\rho}\right)^n > 0 \qquad n>0
\end{equation*}
with $\rho$ and $1-\rho$ paired together. 
\citeauthor{Bombieri} \cite{Bombieri} showed that Li's criterion is an instance of Weil's criterion restricted to a set of test functions $\{g_n(t)\}$ defined as follows.

The Laguerre polynomial  $L_n(t)$ of degree $n\ge0$ and its Laplace transform  are given by
\begin{align*}
 L_n(t) &= \sum_{k=0}^n \frac{1}{k!}\binom{n}{k} (-t)^k \qquad  t\ge0 \\
\widetilde L_n(s) &= -\sum_{k=0}^n \binom{n}{k} \left(-\frac{1}{s}\right)^{k+1}
= \frac{1}{s}\sum_{k=0}^n \binom{n}{k} \left(-\frac{1}{s}\right)^k
= \frac{1}{s}\left(1-\frac{1}{s}\right)^n
\end{align*}
The associated Laguerre polynomial $g_n(t)\equiv L^1_{n-1}(t)$, $\sfg_n(t)\equiv -e^t L^1_{n-1}(-t)$ ($n>0$) and corresponding Laplace transforms are defined by\footnote{Also see \cite{Coffey2, Lagarias} for identification of the test functions as the associated Laguerre polynomials $L^1_{n-1}(t)$.}
\begin{alignat*}{4}
 g_n(t) &= L^1_{n-1}(t)&  &= \sum_{k=0}^{n-1} L_k(t) \\
\wtg_n(s) &= \widetilde{L^1}_{n-1}(s) & &= \frac{1}{s}\sum_{k=0}^{n-1}\left(1-\frac{1}{s}\right)^k
& &= 1-\left(1-\frac{1}{s}\right)^n & &= 1-\left(1-\frac{1}{1-s}\right)^{-n} \\
\wsg_n(s) &= \wtg_n(1-s) & &= \widetilde{L^1}_{n-1}(1-s)
& &=  1-\left(1-\frac{1}{1-s}\right)^n  & &=  1-\left(1-\frac{1}{s}\right)^{-n}
\end{alignat*}
Then, as observed in \cite{Bombieri}, $\wtg_n(s)\wtg_n(1-s) = \wtg_n(s)+\wtg_n(1-s)$ by the identity $(1-r)(1-r^{-1}) \equiv  (1-r)+(1-r^{-1})$.
Hence, as discussed above,  Li' criterion is an instance of Weil's criterion
\begin{equation*}
\lambda_n = \sum_\rho \wtg_n(\rho) = \frac{1}{2}\sum_\rho \wtg_n(\rho)+\wtg_n(1-\rho)
 = \frac{1}{2}\sum_\rho \wtg_n(\rho) \wtg_n(1-\rho) > 0
\end{equation*}
To derive an explicit form,  we use the expressions obtained in the previous section 
\begin{alignat}{1}
  \lambda_n &= \sum_\rho 1-\left(1-\frac{1}{\rho}\right)^n  
 = -\sum_{k=1}^n \binom{n}{k} \sum_\rho\left(-\frac{1}{\rho}\right)^k
 = n\sum_\rho\frac{1}{\rho} - \sum_{k=2}^n \binom{n}{k} \sum_\rho\left(-\frac{1}{\rho}\right)^k
\nonumber \\
 &= n + \frac{n}{2}\left(\gamma - \log4\pi\right)
-\sum_{k=2}^n\binom{n}{k}\left(\eta_{k-1} - (-)^k(1-2^{-k})\zeta(k)\right) - \sum_{k=2}^n\binom{n}{k} (-)^k  
\nonumber \\
\intertext{which, with the aid of the identity 
$\sum_{k=0}^n\binom{n}{k} (-)^k = 1-n+\sum_{k=2}^n\binom{n}{k} (-)^k = 0$, leads to}
 \lambda_n &= 1 + \frac{n}{2}\left(\gamma - \log4\pi\right)
+ \sum_{k=2}^n\binom{n}{k}\left((-)^k(1-2^{-k})\zeta(k) - \eta_{k-1}\right)   
\label{eq:Lieta} \\
\intertext{The corresponding exercise for the expressions involving $\{\mu_k\}$ leads to}
\lambda_n 
&= 1+ \frac{n}{2}\left(\gamma - \log4\pi\right)
+ \sum_{k=2}^n\binom{n}{k}\left( 2^{-k}\zeta(k) + \mu_{k-1}\right)
\label{eq:Limu}
\end{alignat}
\citeauthor{Bombieri} \cite{Bombieri} proved (\ref{eq:Lieta}).
It is an arithmetic interpretation of $\lambda_n$, involving as it does the zeta function at integer argument and  $\{\eta_k\}$ given by the arithmetic formula 
(\ref{eq:eta})\footnote{It is shown in \cite{Bombieri, Coffey} that the sequence $\{\eta_k\}$  can be generated recursively  
with the aid of the Stieltjes constants
\begin{equation*}
 \gamma_k = \frac{(-)^k}{k!} \lim_{x\rightarrow \infty} \left(\sum_{n\le x}\frac{(\log n)^k}{n}   
 - \frac{(\log x)^{k+1}}{k+1} \right)
\end{equation*}
bypassing, thereby, explicit dependence on the von Mangoldt function $\Lambda(n)$.}.

An arithmetic interpretation for the (formally equivalent) form (\ref{eq:Limu}) is less direct because
the sequence $\{\mu_k\}$ (\ref{eq:mu}) arises from analytic continuation rather than from an explicit arithmetic construct.
\citeauthor{Apostol} \cite{Apostol} derived closed form expressions for $\zeta^{(k)}(0)$, from which closed form expressions for $\mu_k$ can be inferred\footnote{Inevitably, $\{\zeta^{(k)}(0)\}$ and hence $\{\mu_k\}$  can also be written in terms of the Stieltjes constants.}.

A study of the properties of $\{\mu_k\}$ might shed light on whether $\lambda_n>0$ for $n>0$ does indeed hold.
If $f(t)$ were a probability density then $\mu_k$ would be the $k^{\rm th}$ moment and therefore necessarily positive.
This does not, of course, hold in this case because $\wtf(s)$ is the Laplace transform of the density $f(t)$ only for $s>1$.
Nonetheless, a probabilistic interpretation arguably seems no less compelling than an arithmetic interpretation.

\section{Discussion}
It may seem natural to ask whether it makes sense for the analytic version (\ref{eq:density2}) of the density $f(t)$ to have an atom at zero when it is meant to be equivalent to the arithmetic version
(\ref{eq:density1}), which is zero at $t=0$ since $\Lambda(1)=0$ (1 is not prime).

In the conditioning notation of probability theory, it is more informative to write the arithmetic density as, say, $f(t|q)$ -- to be understood as
``$f(t)$ given the prime numbers as prior knowledge".
Hence $f(t|q)$ is discrete with atoms at the primes $\log q$, where  $\log 1 = 0$ is excluded by prior construction.
Similarly, the analytic density may be written as $f(t|\rho)$: 
``$f(t)$ given the complex zeros $\rho$ of the zeta function as prior knowledge''.
Such prior knowledge is defined in $s$ space and it forces neither discreteness of $f(t|\rho$) 
nor zero mass at $t=0$.


While $f(t|q)$ and $f(t|\rho)$ are constructed in accordance with their associated prior assumptions, all integrals derived therefrom must necessarily be in agreement. 
It is, after all, the integrals such as the Chebyshev counting function that are the objects of ultimate interest.

A deeper probabilistic approach to the Riemann Hypothesis, drawing specifically from the theory of infinitely divisible distributions, will be explored in a sequel to this paper.
The work discussed here will be of direct relevance because the  Laplace transform plays a central role in the treatment of infinitely divisible distributions on $[0,\infty)$ as discussed in \cite{Feller, SSJS97}.

\section{Conclusion}
We have constructed arithmetic and analytic forms for the density on $[0,\infty)$ from which sums and integrals of interest may be surmised.
We have shown that the general explicit formula follows effortlessly from the specification of such a density.
We have further given an alternative form for Li's criterion.

\bibliography{./XplctForm}

\begin{thebibliography}{10}
\providecommand{\natexlab}[1]{#1}
\providecommand{\url}[1]{\texttt{#1}}
\expandafter\ifx\csname urlstyle\endcsname\relax
  \providecommand{\doi}[1]{doi: #1}\else
  \providecommand{\doi}{doi: \begingroup \urlstyle{rm}\Url}\fi

\bibitem[Abramowitz and Stegun(1965)]{AbrSteg}
M.~Abramowitz and I.~Stegun.
\newblock \emph{Handbook of Mathematical Functions}.
\newblock Dover, New York, 1965.

\bibitem[Apostol(1985)]{Apostol}
T.~M. Apostol.
\newblock Formulas for higher derivatives of the {R}iemann zeta function.
\newblock \emph{Math. Comput.}, 44:\penalty0 223--232, 1985.

\bibitem[Bombieri and Lagarias(1999)]{Bombieri}
E.~Bombieri and J.~C. Lagarias.
\newblock Complements to {L}i's criterion for the {R}iemann {H}ypothesis.
\newblock \emph{J. Number Theory}, 77:\penalty0 274--287, 1999.

\bibitem[Coffey(2004)]{Coffey}
M.~W. Coffey.
\newblock Relations and positivity results for derivatives of the {R}iemann
  $\xi$ function.
\newblock \emph{J. Comput. Appl. Math.}, 166:\penalty0 525--534, 2004.

\bibitem[Coffey(2007)]{Coffey2}
M.~W. Coffey.
\newblock The theta-{L}aguerre calculus formulation of the {L}i/{K}eiper
  constants.
\newblock \emph{J. Approx. Theory}, 146:\penalty0 267--275, 2007.

\bibitem[Edwards(2001)]{Edwards}
H.~M. Edwards.
\newblock \emph{Riemann's Zeta Function}.
\newblock Dover, New York, 2001.

\bibitem[Feller(1971)]{Feller}
W.~Feller.
\newblock \emph{An Introduction to Probability Theory and its Applications,
  Vol. II}.
\newblock Wiley, New York, 1971.

\bibitem[Lagarias(2007)]{Lagarias}
J.~C. Lagarias.
\newblock Li coefficients for automorphic {L}-functions.
\newblock \emph{Ann. Inst. Fourier}, 57:\penalty0 1689--1740, 2007.

\bibitem[Li(1997)]{Li}
X.-J. Li.
\newblock The positivity of a sequence of numbers and the {R}iemann
  {H}ypothesis.
\newblock \emph{J. Number Theory}, 65:\penalty0 325--333, 1997.

\bibitem[Sibisi and Skilling(1997)]{SSJS97}
S.~Sibisi and J.~Skilling.
\newblock Prior distributions on measure space.
\newblock \emph{J. R. Statist. Soc. B}, 59:\penalty0 217--235, 1997.

\end{thebibliography}
\bibliographystyle{plainnat}
\end{document}